\documentclass[12pt,a4paper,draft]{article}
\usepackage[english,russian]{babel}
\usepackage{amssymb, amsmath, latexsym, amsfonts, amsthm}
\begin{document}
\renewcommand{\refname}{References}
\newtheorem{theorem}{Theorem}
\newtheorem{lemma}{Lemma}
\newtheorem{corollary}{Corollary}
\newtheorem{proposition}{Proposition}
\begin{center}
On a class of infinitely differentiable functions in ${\mathbb R}^n$ 
admitting holomorphic extension in ${\mathbb C}^n$ 
\end{center}
\begin{center}
I.Kh. Musin\footnote {supported by the Russian Science Foundation grant 15-01-01661 and RAS Program of Fundamental Research "`Modern problems of theoretical mathematics"'}, P.V. Yakovleva\footnote {supported by the Russian Science Foundation grant 15-01-01661}
\end{center}

\renewcommand{\abstractname}{}
\begin{abstract}
{\sc Abstract}. A space $G(M, \varPhi)$ of infinitely differentiable functions in ${\mathbb R}^n$ constructed with a help of a family $\varPhi=\{\varphi_m\}_{m=1}^{\infty}$ of real-valued functions 
$\varphi_m \in~C({\mathbb R}^n)$ and a logarithmically convex sequence $M$ of positive numbers is considered in the article. In view of conditions on $M$ each function of $G(M, \varPhi)$ can be extended to an entire function in ${\mathbb C}^n$. Imposed conditions on $M$ and $\varPhi$ allow to describe the space of such extensions. 

\vspace {0.3cm}






\vspace {0.3cm}
\end{abstract}

\vspace {0.6cm}

{\bf 1. On the problem}. L. Ehrenpreis \cite {E}, B.A. Taylor \cite {T} considered spaces of infinitely differentiable functions in ${\mathbb R}^n$ constructed by the following general scheme. 
Let $M=(M_k)_{k=0}^{\infty}$ be a sequence of positive numbers $M_k$ such that $M_0=1$ and 
$
\displaystyle\lim_{k\rightarrow\infty}\frac {\ln M_k} {k} = +\infty.
$
Let $\varPhi=\{\varphi_m\}_{m=1}^{\infty}$ be a family of real-valued functions 
$\varphi_m \in~C({\mathbb R}^n)$ such that for each $m \in {\mathbb N}$:

1). $
\lim \limits_{x \to \infty}
\displaystyle \frac {\varphi_m(x)} {\Vert x \Vert} = +\infty $
($\Vert \cdot \Vert$ be the Euclidean norm in ${\mathbb R}^n$);

2). $\lim \limits_{x \to \infty}
(\varphi_m(x) - \varphi_{m+1}(x)) = +\infty $.

\noindent Define a space $G(M, \varPhi)$ of infinitely differentiable functions $f$ in ${\mathbb R}^n$ as follows. For each $m\in\mathbb N$ and $\varepsilon > 0$ let $G_{m,\varepsilon}$ be a space  of $C^{\infty}$-functions $f$ in ${\mathbb R}^n$ with a finite norm
$$ 
p_{m, \varepsilon}(f)= \sup_{x \in {\mathbb R}^n, \alpha \in {\mathbb Z_+^n}} \frac {\vert (D^{\alpha}f)(x)\vert}{\varepsilon^{\vert \alpha \vert} M_{\vert \alpha \vert} e^{\varphi_m(x)}} \ .
$$ 
Put $G(M, \varPhi)=\displaystyle \bigcap_{m=1}^{\infty}\bigcap_{\varepsilon > 0} G_{m,\varepsilon}$. 
With usual operations of addition and multiplication by complex numbers $G(M, \varPhi)$ is a linear space. 
The family of norms $p_{m, \varepsilon}$ defines a locally convex topology in $G(M, \varPhi)$.

Under additional conditions on $\varPhi$ and $M$ various problems of harmonic analysis and operator theory in $G(M, \varPhi)$ were considered in \cite {E}, \cite {T}, \cite {M}. 

Here we study the problem of holomorphic extension of functions belonging to $G(M, \varPhi)$. 
Namely, our aim is to impose conditions on $\varPhi$ and $M$ which guarantee that all functions belonging to $G(M, \varPhi)$ admit extension to entire functions and allow to describe the space of these extensions. Similar problems  we considered earlier in \cite {MY} in case of a subspace of a space of rapidly decreasing smooth functions on an unbounded closed convex set in ${\mathbb R}^n$.

Note that in \cite {A} there was considered the problem of extension of infinitely differentiable functions in ${\mathbb R}^n$ satisfying with all their partial derivatives some weighted estimates to entire functions with consistent estimates for growth. Moreover, weight functions in \cite {A} depend on modules of variables. 

{\bf 2. Conditions on $\varPhi$ and $M$ and the main result.}
We assume that
the sequence  $M$ satisfies the following additional conditions:

$\alpha_1)$. $M_k^2 \le M_{k-1}M_{k+1},  \ \ \forall k \in {\mathbb N}$; 


$\alpha_2)$. $ \forall \varepsilon > 0$ $ \exists a_{\varepsilon} > 0 $ $\forall k \in {\mathbb Z_+}$ \
$M_k \le a_{\varepsilon}\varepsilon^k k!$;

$\alpha_3)$. there exists a logarithmically convex sequence 
$K=(K_m)_{m=0}^{\infty}$ with $K_0=1$ such that for some $t_1>1, t_2 > 1$ 
$$
t_1^{-1} t_2^{-m} K_m \le \frac {m!}{M_m} \le t_1 t_2^m K_m, \ \forall m \in {\mathbb Z_+}.
$$

Denote by $w_K$ a function defined as follows:
$w_K(r) = \displaystyle \sup_{m \in {\mathbb Z_+}} \ln\frac {r^m} {K_m}$ for $r>0$, $w_K(0) = 0$. 

We impose the following condition on $\varPhi$: for each $m \in {\mathbb N}$ and for each $\delta >0$ there exist  constants $a_m > 1$ (depending only on $m$) and $b_{m, \delta} > 0$ such that for each $R>0, x \in {\mathbb R}^n$ and
$\xi \in {\mathbb R}^n$ with $\Vert \xi \Vert \le R$
\begin{equation}
\varphi_{m + 1}(x + \xi) \le \varphi_m (x) + \omega_K(a_m \delta R) + b_{m, \delta}.
\end{equation}

Let $E(K, \varPhi)$ be the space of entire functions $f$ in ${\mathbb C}^n$ such that for each $m \in {\mathbb N}$ and for each $\varepsilon > 0$ there exists a constant $c_{m, \varepsilon} > 0$ such that
$$
\vert f(z) \vert \le c_{m, \varepsilon} e^{\varphi_m(x) + w_K(\varepsilon \Vert y \Vert)}, \ z \in {\mathbb C}^n,
$$
endowed with a topology defined by the family of norms 
$$
q_{m, \varepsilon}(f) = \sup_{z \in {\mathbb C}^n} 
\frac {\vert f(z)}{e^{\varphi_m(x) + w_K(\varepsilon \Vert y \Vert)}} \ , \ \varepsilon > 0, m \in {\mathbb N}.
$$
Here as usual $x = Re  z, y = Im  z$. 


Using elementary methods we prove the following theorem. 

{\bf Theorem.} {\it The spaces $G(M, \varPhi)$ and $E(K, \varPhi)$ are isomorphic}.

\vspace {0.01cm}

{\bf 3. Auxiliary result}. The following Lemma will be used in the proof of the Theorem.

{\bf Lemma.} 
{\it For each $r \ge 0$ \ 
$$
2w_K(r) \le w_K(e t_2^2 r) + 3 \ln t_1. 
$$}

{\bf Proof}. Note first that in view of logarithmical convexity of $M$ for all  $p, q \in {\mathbb Z_+}$ 
\  $M_{p+q} \ge M_p M_q$. So and also taking into account the condition $\alpha_3)$ 
we have for all $p, q \in {\mathbb Z_+}$
$$
K_{p+q} \le t_1 t_2^{p+q} \frac {(p+q )!}{M_{p+q}} \le t_1 (t_2 e)^{p+q} \frac {p! q!}{M_pM_q} \ .
$$
Again using the condition $\alpha_3)$ we get
$$
K_{p+q} \le t_1^3 (et_2^2)^{p+q}K_p K_q  \ . 
$$
With help of this inequality we have for $r>0$ 
$$
2w_K(r) = 
\sup_{m \in {\mathbb Z_+}} \ln \frac {r^{2m}}{K_m^2} \le
\sup_{m \in {\mathbb Z_+}} \ln \frac {t_1^3 (e t_2^2 r)^{2m}}{K_{2m}} 
\le w_K(e t_2^2 r) + 3 \ln t_1.
$$
For $r=0$ the assertion of lemma is obvious. So lemma is proved. 

{\bf 4. Proof of the Theorem}.  
Let $f~\in~G(M, \varPhi)$. Then 
$\forall m \in {\mathbb N}$, $\varepsilon > 0$
\begin{equation}
\vert (D^{\alpha}f)(x)\vert \le p_{m, \varepsilon}(f) 
\varepsilon^{\vert \alpha \vert} M_{\vert \alpha \vert} e^{\varphi_m(x)}  \ , \ x \in {\mathbb R}^n, \alpha \in {\mathbb Z_+^n}.  
\end{equation}
From this and the condition $\alpha_2)$ it follows that 
$$
f(x)=\displaystyle \sum_{\vert \alpha \vert \ge 0} \frac {(D^{\alpha}f)(x_0)}{\alpha!} (x-x_0)^{\alpha}, \ x, x_0 \in {\mathbb R}^n,
$$
and the series on the right hand side of this equality converges uniformly on compact subsets of ${\mathbb R}^n$ to $f$.
In view  of (2) and the condition $\alpha_2)$ on $M$ for each $x_0 \in {\mathbb R}^n$ the function 
$F_{f, x_0}$ defined as 
$$
F_{f, x_0}(z)= \displaystyle \sum_{\vert \alpha \vert \ge 0} \frac {(D^{\alpha}f)(x_0)}{\alpha!} (z-x_0)^{\alpha} \ , \ z \in {\mathbb C}^n,
$$
is an entire function.
Note that for $x \in {\mathbb R}^n$ \ $F_{f, x_0}(x)=f(x)$. So for $x_1, x_2 \in {\mathbb R}^n$ \ 
$F_{f, x_1}(z)= F_{f, x_2}(z), \ z \in {\mathbb C}^n$. 
Thus, for each $f \in G(M, \varPhi)$ we defined a function $F_f \in H({\mathbb C}^n)$ such that for each $\xi \in {\mathbb R}^n$ we have  $F_f= F_{f, \xi}$ in ${\mathbb C}^n$ and $F_f(x)=f(x), \ x \in {\mathbb R}^n$.  

To show that  $F_f \in E(K, \varPhi)$ we need to estimate a growth of $F_f$.
Since 
$$
F_f(z) = \displaystyle \sum_{\vert \alpha \vert \ge 0} \frac 
{(D^{\alpha}f)(x)}{\alpha!} (iy)^{\alpha} , \ z = x + iy, x, y \in {\mathbb R}^n,
$$
then for each $m \in {\mathbb N}$ and $\varepsilon > 0$
$$ 
\vert F_f(z) \vert \le p_{m, \varepsilon}(f)
\sum_{\vert \alpha \vert \ge 0} \frac 
{\varepsilon^{\vert \alpha \vert} M_{\vert \alpha \vert} \Vert y \Vert^{\alpha} e^{\varphi_m(x)}} {\alpha!} = 
$$
$$
= 
p_{m, \varepsilon}(f) e^{\varphi_m(x)}
\sum_{N=0}^{\infty}
\varepsilon^N M_N \Vert y \Vert^N 
\sum_{\vert \alpha \vert = N} \frac {1}{\alpha!} = 
$$
$$
= p_{m, \varepsilon}(f) e^{\varphi_m(x)}
\sum_{N=0}^{\infty}
\varepsilon^N M_N \Vert y \Vert^N \frac {n^N}{N!} \le 
2 t_1 p_{m, \varepsilon}(f) e^{\varphi_m(x)}
\sup_{N \in {\mathbb Z_+}}
\frac 
{(2\varepsilon n t_2 \Vert y \Vert)^N} {K_N} =  
$$
$$
=2 t_1 p_{m, \varepsilon}(f) e^{\varphi_m(x) + w_K(2\varepsilon n t_2 \Vert y \Vert)}. 
$$

Define a map $T: G(M, \varPhi) \rightarrow E(K, \varPhi)$ by the formula $T(f)=F_f, f \in G(M, \varPhi)$. 
Obviously, the mapping $T$ is one-to-one and linear. 
From the last inequality it follows that for all $\varepsilon > 0, m\in\mathbb N$
$$
q_{m, 2\varepsilon n t_2}(T(f)) \le 2 t_1 p_{m,\varepsilon}(f) \ , \ f \in G(M, \varPhi).
$$
This means that  $T$ is a continuous mapping from $G(M, \varPhi)$ to $E(K, \varPhi)$. 

Now we prove that $T$ is surjective and the inverse mapping $T^{-1}$ is continuous. 
For that take $F \in E(K, \varPhi)$ and show that $f=F_{|{\mathbb R}^n}$ belongs to $G(M, \varPhi)$. 
By the integral Cauchy formula for each $m \in {\mathbb N}$, $R>0$, $\alpha \in {\mathbb Z_+^n}$ we have that 
$$
(D^{\alpha}f)(x) = 
\frac {\alpha! }{(2\pi i)^n} 
\displaystyle \int \cdots \int_{L_R(x)}
\frac 
{F(\zeta)}
{(\zeta_1 - x_1)^{\alpha_1 +1} \cdots (\zeta_n - x_n)^{\alpha_n +1}} \ d \zeta , \ x \in {\mathbb R}^n,
$$
where  
$L_R(x)= \{\zeta = (\zeta_1, \ldots , \zeta_n) \in {\mathbb C}^n: \vert \zeta_j - x_j \vert = R, j=1, \ldots , n \}$, 
$d \zeta = d \zeta_1 \cdots d  \zeta_n $.
From this we get that
$$
\vert (D^{\alpha} f)(x) \vert 
\le \frac 
{\alpha!}{R^{\vert \alpha \vert}}  
\max_{\zeta \in L_R(x)} \vert F(\zeta) \vert.
$$ 
Obviously, 
$$
 \vert (D^{\alpha} f)(x) \vert 
\le  \frac 
{\alpha! q_{m + 1, \varepsilon}(F)}{R^{\vert \alpha \vert}}  
\exp ({\max \limits_{\Vert \xi \Vert \le \sqrt n R} \varphi_{m + 1}(x + \xi) + w_K(\varepsilon \sqrt n R))}.
$$ 
By the additional condition on $\varPhi$ (inequality (1)) we have that 
$$
\vert (D^{\alpha} f)(x) \vert 
\le  \frac 
{\alpha!q_{m + 1, \varepsilon}(F)}{R^{\vert \alpha \vert}}  
e^{\varphi_m(x) + 2 \omega_K(\varepsilon a_m \sqrt n R) + b_{m, \varepsilon}}.
$$ 
Using Lemma we can find positive constants $c_m$ and $d_{m, \varepsilon}$ such that 
$$
\vert (D^{\alpha} f)(x) \vert 
\le  \frac 
{\alpha!q_{m + 1, \varepsilon}(F)}{R^{\vert \alpha \vert}}  
e^{\varphi_m(x) + \omega_K(c_m \varepsilon R) + d_{m, \varepsilon}}.
$$ 
From this we have that
$$
\vert (D^{\alpha} f)(x) \vert 
\le \alpha! q_{m + 1, \varepsilon}(F) 
\inf \limits_{R > 0}
\frac {e^{\omega_K(c_m \varepsilon R)}}
{R^{\vert \alpha \vert}}  
e^{\varphi_m(x) + d_{m, \varepsilon}}.
$$ 
Hence,
$$
\vert (D^{\alpha} f)(x) \vert 
\le \alpha! q_{m + 1, \varepsilon}(F) 
(c_m \varepsilon)^{\vert \alpha \vert}
\inf \limits_{r > 0}
\frac {e^{\omega_K(r)}}
{r^{\vert \alpha \vert}}  
e^{\varphi_m(x) + d_{m, \varepsilon}}.
$$ 
From this we have that
$$
\vert (D^{\alpha} f)(x) \vert 
\le \frac {\vert \alpha \vert!}{K_{\vert \alpha \vert}} q_{m + 1, \varepsilon}(F) (c_m \varepsilon)^{\vert \alpha \vert}
e^{\varphi_m(x) + d_{m, \varepsilon}}.
$$ 
Consequently,
$$
\vert (D^{\alpha} f)(x) \vert 
\le t_1 e^{d_{m, \varepsilon}} q_{m + 1, \varepsilon}(F) 
(c_m t_2 \varepsilon)^{\vert \alpha \vert} M_{\vert \alpha \vert}  
e^{\varphi_m(x)}.
$$ 
Hence,
$$
p_{m, c_m t_2 \varepsilon}(f) \le t_1 e^{d_{m, \varepsilon}} q_{m + 1, \varepsilon}(F).
$$
This means that $f=F_{|{\mathbb R}^n} \in G(M, \varPhi)$ and the mapping $T^{-1}$ is continuous. Clearly, $T(f)=F$. Thus, the spaces $G(M, \varPhi)$ and $E(K, \varPhi)$ are isomorphic.~$\square$

\vspace {0.3cm}

Il'dar Kh. Musin 

Institute of Mathematics with Computer Centre of Ufa Scientific Centre of Russian Academy of Sciences; 
Chernyshevsky str., 112, Ufa, 450077, Russia

musin\_ildar@mail.ru

\vspace {0.3cm}

Polina V. Yakovleva

Ufa State Aviation Technical University; 
K. Marx str., 12, Ufa, 450000, Russia 

polina81@rambler.ru

\end{document}